\documentclass{article}

\usepackage{amssymb,amsmath,amsthm}
\usepackage{amsfonts}
\usepackage{enumerate}
\usepackage[T1]{fontenc}
\usepackage[sort,numbers]{natbib}
\def\car{{\mathbf 1}}
\def\N{{\mathbf N}}
\def\R{{\mathbf R}}
\def\P{{\mathbf P}}
\def\d{\text{d}}
\def\esp#1{{\mathbf E}\left[#1\right]}
\newtheorem{defn}{Definition}
\newtheorem{thm}{Theorem}
\begin{document}

  \title{Upper bound of loss probability  in an OFDMA system with
    randomly located users}
 \author{L. Decreusefond and E. Ferraz and     P. Martins\thanks{All
     authors are with ENST/ CNRS UMR 5141, 46, rue Barrault, Paris, FRANCE}}

  \maketitle

\begin{abstract}
  For OFDMA systems, we find a rough but easily computed upper bound for the probability
  of loosing communications by insufficient number of sub-channels on
  downlink. We consider as random the positions of receiving users in the system as
  well as the number of sub-channels dedicated to each one. We use
  recent results of the theory of point processes which reduce our
  calculations to that of the first and second moments of the total
  required number of sub-carriers. 
\end{abstract}

\section{Introduction}

The demand for high data rate wireless applications with restrictions
in the RF signal bandwidth requires bandwidth efficient air interface
schemes. It is known that OFDM yields a relatively simple solution to
these problems~\cite{stuber}. Based on the OFDM system, OFDMA can
achieve a larger capacity. Furthermore, this latter system is more flexible,
since it can be easily scaled to fit in a certain piece of spectrum
simply by changing the number of used
subcarriers~\cite{vannee}. However, as any wireless systems, OFDM and
OFDMA have physical limitations which cause loss of communications. This
loss can be caused by insufficient power  or by low
signal-to-interference ratio, for instance. In this paper we are
interested in the calculation of  an upper bound of the probability of loosing
a communication due to an insufficient number of sub-channels in the
downlink.

We say that the system is overloaded when all non-used sub-channels
are not enough to warrant a minimum data rate for an incoming
demand. We consider a system with $N_0$ sub-carriers and $N_i$ is the number of sub-carriers used
by the $i$-th user in the cell. As it is usually done, we substitute
the finite number of subcarriers by  infinity and substitute the loss
probability  by $P_{loss}=P(\sum_i N_i>N_0)$. It is well known
that this consideration gives us an upper bound for the actual loss
probability.

A user $i$ requires a capacity $C_i$ depending on the service he uses.
Considering a system with just one kind of service, all users require
the same capacity $C_0$. Even so, the number of subcarriers for each
user varies according to the channel conditions. These conditions can
be summarized into two kinds of gains,  one depending
only on the position of the user $i$, the path loss $G_{pl_i}$, and a
 gain $G_i$, which classically may include the shadowing and the
Rayleigh fading. We choose the simplest model to represent the path
loss:
\begin{equation*}
  G_{pl_i}=\frac{K}{D_i^{\gamma}}
\end{equation*}
where $K$ and $\gamma$ are constants and $D_i$ is the distance between
the user $i$ and the antenna.

Shannon's maximum achievable capacity implies that:
\begin{eqnarray*}
  N_i=\left\lceil\frac{C_0}{W\log_2\left(1+\dfrac{P_t G_{pl_i}G_i}{I}\right)}\right
  \rceil
\end{eqnarray*}
where $W$ is the bandwidth of each sub-carrier, $P_t$ is the mean
transmitted power by sub-channel and $I$ is the power of the additive
Gaussian white noise by sub-channel.

We consider that the number as well as the position of users in the cells are random. After some natural assumptions done in the
following section, we conclude that the configuration of users in the
cell is a Poisson point process (see section~\ref{ppp}).

After a summary on Poisson point process, we consider three different
cases to calculate an upper bound for the loss probability. First we
consider the simplest case with deterministic
gain. In Section~\ref{random one}, we  consider a non-selective frequency
gain, the shadowing. In section~\ref{general}, we consider a general case
from which all other cases could be derived but for which no closed
form formula exists.

\section{Poisson point processes}
\label{ppp}

For details on point processes, we refer to  
\cite{MR1950431,kallenberg83,kallenberg,MR1113698}. A configuration $\eta$ in $\R^k$ is a set $\{x_n,\, n\ge 1\}$ where
for each $n\ge 1$, $x_n\in \R^k$, $x_n\neq x_m$ for $n\neq m$ and each compact subset of
$\R^k$ contains only a finite subset of $\eta$. We denote by
$\Gamma_{\R^k}$ the set of configurations in $\R^k$. Equipped with the
vague topology of discrete measures, $\Gamma_{\R^k}$ is a complete,
separable metric space. A point process
$\Phi$ is a random variable with values in $\Gamma_{\R^k}$, i.e.,
$\Phi(\omega)=\{X_n(\omega),\, n\ge 1\}\in \Gamma_{\R^k}$. For
$A\subset \R^k$, we denote by $\Phi_A$ the random variable which
counts the number of atoms of $\Phi(\omega)$ in $A$:
\begin{equation*}
  \Phi_A(\omega)=\sum_{n\ge 1}\car_{X_n(\omega)\in A}\ \in \N\cup \{+\infty\}.
\end{equation*}
Poisson point processes are particular instances of point processes
such that:
\begin{defn}
\label{def_poisson_process}
  Let $\Lambda$ be a $\sigma$ finite measure on $\R^k$. A point
  process $\Phi$ is a Poisson process of intensity $\Lambda$ whenever
  the two following properties hold.
  \begin{enumerate}[1 - ]
  \item For any compact subset $A\in \R^k$, $\Phi_A$ is a Poisson random
    variable of parameter $\Lambda(A)$, i.e.,
    \begin{equation*}
      \P(\Phi_A=k)=e^{-\Lambda(A)}\frac{\Lambda(A)^k }{k!}.
    \end{equation*}
\item For any disjoint subsets $A$ and $B$, the random variables
  $\Phi_A$ and $\Phi_B$ are independent.
  \end{enumerate}
\end{defn}

The notion of point process trivially extends to configurations in
$\R^k\times X$ where $X$ is a subset of $\R^m$. A configuration is
then typically of the form $\{(x_n,\, y_n),\, n\ge 1\}$ where
for each $n\ge 1$, $x_n\in \R^k$ and $y_n\in X$. We keep writing $(x_n,\, y_n)$ as a
couple, though it could be thought as an element of $\R^{k+m}$, to stress the asymmetry between the spatial coordinate $x_n$
and the so-called mark, $y_n$. For a marked point process, we denote
by $\Phi$ the set of locations, i.e., $\Phi(\omega)=\{X_n, \, n\ge
1\}$ and by $\bar{\Phi}$ the set of both locations and marks,  i.e.,
$\bar{\Phi}(\omega)=\{(X_n,\, Y_n), \, n\ge
1\}$. A marked point process with position
dependent marking is a marked point process for which the law of
$Y_n$, the mark associated to the atom located at $X_n$, depends only
on $X_n$ through a kernel $K$:
\begin{equation*}
  \P(Y_n\in B\, |\, \Phi)=K(X_n,\, B), \text{ for any } B\subset X.
\end{equation*}
If $K$ is a probability kernel, i.e., if $K(x,\, X)=1$ for any $x\in
\R^k$ then it is well known that $\bar{\Phi}$ is a Poisson process of
intensity $K(x,\, \d y)\d\Lambda(x)$ on $\R^k \times \R^m$.
The Campbell formula is a well known and useful formula
\begin{thm}
\label{thm_campbell}
Let $\bar{\Phi}$ be a marked Poisson process on $\R^k\times \R^m
$. Let  $\Lambda$ be the intensity of the underlying Poisson process
and $K$ the kernel of the position dependent marking.
  For  $f\, :\, \R^k\times \R^m\to \R$ a measurable non-negative function, let 
  \begin{equation*}
    F=\int f\, \d\bar\Phi = \sum_{n\ge 1}f(X_n,\, Y_n ).
\end{equation*}
Then,
\begin{equation*}
  \esp{F}=\int_{\R^k\times \R^m} f(x,\, y)\, K(x,\, \d y)\d \Lambda(x).
\end{equation*}
\end{thm}
\begin{defn}
  For $F\, : \, \Gamma_{\R^k}\to \R$, for any $x\in \R^k$, we define
  \begin{equation*}
    D_xF(\omega)=F(\omega\, \cup\, \{x\})-F(\omega).
  \end{equation*}
\end{defn}
Note that for $F=\int f\d\Phi$, $D_xF=f(x),$ for any $x\in \R^k$.
We now quote from \cite{MR1962538,Wu:2000lr}  the main result on which our
inequalities are based:
\begin{thm}[Concentration inequality]
\label{thm_concentration}
Assume that $\Phi$ is a Poisson process on $\R^k$ of intensity $\Lambda.$
  Let $f\, :\, \R^k\to \R^+$ a measurable non-negative function and let 
   \begin{equation*}
    F(\omega)=\int f\, \d\Phi = \sum_{n\ge 1}f(X_n(\omega)).
  \end{equation*}
Assume that $|D_xF(\omega)|\le s$ for any $x\in \R^k$.
Let 
\begin{equation*}
  m_F=\esp{F}=\int f(x)\, \d \Lambda(x) 
\end{equation*}
and 
\begin{equation*}
  v_F= \int |D_xF(\omega)|^2\d \Lambda(x)=\int f^2(x)\, \d \Lambda(x).
\end{equation*}
Then, for any $t\in \R^+$,
\begin{equation*}
  \P(F-m_F\ge t)
\le\\ \exp\left(-\, \frac{v_F}{s^2}\ g\left(\frac{t\, s}{v_F}\right)\right)
\end{equation*}
where $g(t)=(1+t)\ln(1+t)-t$. 
\end{thm}
\section{Deterministic gain}
\label{deterministic}
We state the following assumptions:
\newtheorem{assumption}{Assumption}
\begin{assumption}
  \label{a1}
  The position of each user is  independent
  on the position of all other. The users are indistinguishable, i.e.,
  the positions are identically distributed.
\end{assumption}
\begin{assumption}
  \label{a2}
  The time between two consecutive demands of users for service in the
  system (or interarrival time) is exponentially distributed.
\end{assumption}
We define $\rho(x)$ as the surface density of interarrival time in
s$^{-1}$m$^{-2}$, constant in time. Hence, for a region $H\subseteq
B$, the mean interarrival rate is $h=\int_H\rho(x)dx$ in s$^{-1}$.
\begin{assumption}
  \label{a3}
  The service time for every user is exponentially distributed with
  mean $1/\nu$.
\end{assumption}
\begin{assumption}
  \label{a4}
  The cell $C$ is circular, with radius $R$ and with the antenna in
  the center.
\end{assumption}
\begin{assumption}
  \label{a5}
  The channel gain depends only on the distance from the transmitting
  antenna.
\end{assumption}
\begin{assumption}
  \label{a6}
  The surface density of interarrival time is constant.
\end{assumption}
These assumptions are commonly done to simplify the mathematical
treatment and are quite reasonable. 
If we show that the point process given by the location of the users
is a Poisson process, then it is sufficient to  have the two first
moments in order to apply theorem  \ref{thm_concentration} and then  calculate an upper bound $P_{sup}$ for
$P_{loss}$. To do this, we consider the following lemma:
\newtheorem{lemma}{Lemma}
\begin{lemma}
  \label{lemmappp}
  Considering assumptions~\ref{a1}, \ref{a2} and~\ref{a3}, the point
  process  $\Phi$ of the active users positions is, in equilibrium, a
  Poisson process with intensity $\d\Lambda(x)=\rho(x)\nu^{-1}\d x$
\end{lemma}
\begin{proof}
For a region $H$, in virtue of assumptions 2 and 3, the number of receiving (i.e.,
active) customers is the same as the number of customers in an M/M/$\infty$
queue with input rate $h$ and mean service time $\nu^{-1}.$
It is known \cite{MR1996883} that  the distribution of the
  number of users $U$ in equilibrium is then
  \begin{eqnarray*}
    P(U=u)=\frac{(h/\nu)^u}{u!}e^{-h/\nu}.
  \end{eqnarray*}
It follows that the first  condition  of definition \ref{def_poisson_process} is satisfied with intensity measure
  $\Lambda(H)$  \begin{eqnarray*}
    \Lambda(H)=h/\nu=\int_H \frac{\rho(x)}{\nu}\ \d x.
  \end{eqnarray*}
  Condition 2 of definition \ref{def_poisson_process} follows straightforwardly from
  assumption~\ref{a1}.
\end{proof}
Without loss of generality, we consider the cell $C$ has its antenna
located at the origin. We are looking at evaluating
\begin{equation*}
  \P(\int N \, \d \Phi \ge N_0),
\end{equation*}
where $N(x)$ is defined by 
\begin{equation*}
  N(x)=\left\lceil \frac{C_0}{W\log_2 \left(1+\dfrac{P_t K \bar{g}}{I x^{\gamma}}
      \right)}  \right\rceil,
\end{equation*}
where $\bar{g}$ is the mean gain due to shadowing.
Note that, with respect to $x$, $N$ is increasing and piecewise
constant. Let $R_j,\, j=1,\, \cdots,
N_{max}$ be the values such that $N(x)=j$ for $x\in [R_j,\, R_{j+1}).$ We
can easily determine them by 
\begin{equation*}
  R_j=\left(\frac{P_t K\overline{g}}{I(2^{C_0/(jW)}-1)}\right)^{1/\gamma}.
\end{equation*}
According to Theorem \ref{thm_campbell}, it is then clear that 
\begin{equation*}
  \esp{\int N \, \d \Phi}=\int N\d \Lambda=\frac{\pi\rho}{\nu}\sum_{j=1}^{N_{max}}j(R_j^2-R_{j-1}^2) .
\end{equation*}
We denote by $m_N$ the last quantity. Moreover,
\begin{equation*}
  \int N^2\d \Lambda=\frac{\pi\rho}{\nu}\sum_{j=1}^{N_{max}}j^2(R_j^2-R_{j-1}^2).
\end{equation*}
We denote by $v_N$ the last quantity.
We take $N_0$ of the form $\alpha m_N$, so that according to Theorem \ref{thm_concentration}:
\begin{equation*}
  \P(\int N\, \d \Phi\ge \alpha m_N)\le P_{sup}(\alpha)
\end{equation*}
where
\begin{equation*}
  \label{psup}
P_{sup}(\alpha)=
\exp\left(-\frac{v_N}{N_{max}^2}\ g\left(\frac{(\alpha-1)m_N N_{max}}{v_F}\right)\right).
\end{equation*}
It is then natural to verify how far this bound is from the exact
value of the loss probability in simple situations where simulation is available.
 We used here  $\gamma=2.8$, $C_0=200$ kb/s, $W=250$ kHz
and $P_tK/I=1\times 10^6$.  For the surface density of
interarrival time we use $\rho=0.0006$ min$^{-1}$m$^{-2}$ and the
service time  is $1/\nu=1$ min, so, the mean number of users in
the system is $\pi R^2 \rho/\nu=18.85$ users. If we consider the
shadowing with $\sigma=\sqrt{10}$ dB and $\mu=6$ dB, we can use the
mean gain $g$, giving $\overline{g}=1/12$. Thus, users in the cell
boundary use 3 sub-channels, so $N_{max}=3$. 
For $\alpha$ varying from $1$ to $2$, which corresponds here to loss
probabilities about $2\%$ or $0.01\%$, we computed  $\Delta=\log_{10}
P_{sup}/P_{loss}$.  
\begin{table}
  \centering
  \begin{tabular}{ |c|c|c|c|c|c|c|}
\hline
 $\alpha$  & 1.5 & 1.6 &1.7 & 1.8 & 1.9 & 2\\\hline\hline
$P_{sup}$ &    0.18 &   0.1  &  0.04 &   0.02 & 0.008 & 0.003\\\hline
$\Delta$  & 0.98&  0.1 & 1.15 & 1.3 & 1.3 & 1.4\\
\hline
  \end{tabular}
   \caption{Comparison between $P_{sup}$ and $P_{loss}$ for
     deterministic gain.}
  \label{tab:deterministic}
\end{table}
Though concentration inequalities are usually thought as almost
optimal, the results shown in Table \ref{tab:deterministic} seem at
first glance disappointing. Remind though that the computation of the
bound is immediate whereas the simulation on a fast PC took several
hours to get a decent confidence interval. Remind also that the error
is about the same order of magnitude as the error made when using a
usual trick which consists in replacing infinite buffers by finite
ones in  Jackson networks (see \cite{decreusefond93_6}). The margin provided by the bounds may be
viewed as a protection against errors in the modelling or in the
estimates of the parameters.

\section{Random gain}
\label{random one}
Let us determine now the upper bound probability $P_{sup}$ for $P_{loss}$ 
without assumption~\ref{a5} but holding all other assumptions of the preceding
section. Lemma~\ref{lemmappp} still holds, since it is a consequence of 
assumptions~\ref{a1}, \ref{a2} and~\ref{a3}. We also state two other natural 
assumptions:

\begin{assumption}
\label{a8}
The random gain is totally described by the log-normal shadowing, with mean 
$\mu$ and standard deviation $\sigma$, both in dB. 
\end{assumption}
For a user at distance $d$ from the origin, the gain is $G= 1/S$,
where $S$ follows a log-normal distribution:
\begin{eqnarray*}
p_{S}(y)=\frac{\xi}{\sqrt{2\pi}\sigma y}
\exp\left[-\frac{(10\log_{10}y-\mu)^2}{2\sigma^2}\right],
\end{eqnarray*}
where $\xi=10/\ln 10$. 
\begin{assumption}
\label{a7}
A user is able to receive the signal only if the signal-to-interference ratio
is above some constant $\beta_{min}$. 
\end{assumption}
This means, in particular, that the number of subcarriers needed by a
transmitting user is surely  bounded
by 
\begin{equation*}
  N_{max}=\left\lceil \frac{C_0}{W\log_2 (1+\beta_{min})}
      \right\rceil.
\end{equation*}
The situation is slightly different from that of Section
\ref{deterministic}, since the
functional depends on two  aleas: positions and
gains. Consider now that our configurations are of the form $(x,s)$
where $x\in \R^2$ is still a position and $s\in \R$ is a gain. Since
gain and positions are independent, we then have a Poisson process on
$\R^{3}$ of intensity measure $ d\Lambda(x) \otimes p_S(y)\d
y$. Thus we want to evaluate an upper bound of 
\begin{equation*}
  \P(\int N \d \Phi \ge \N_0)
\end{equation*}
where 
\begin{equation*}
  N(x,\, y)=\left\lceil \frac{C_0}{W\log_2 \left(1+\dfrac{P_t K }{I y x^{\gamma}}
      \right)}  \right\rceil.
\end{equation*}
According to Theorem \ref{thm_concentration}, we must  compute 
 \begin{equation*}
   m_N=\int N(x,\ y) p_S(y)\d
y \ \d \Lambda(x)
 \end{equation*}
and
\begin{multline*}
  v_N=\sup_\omega \int |D_{x,y}F(\omega)|^2 p_S(y)\d
y \ \d \Lambda(x)\\
=\int N^2(x,\ y) p_S(y)\d
y \ \d \Lambda(x).
\end{multline*}
Let $\beta_0=\infty$ and $\beta_j=2^{C_0/(Wj)}-1$ for $j=1, \cdots,\,
N_{max}-1$. For $j=1,\cdots,\, N_{max}-1$, let
\begin{equation*}
  A_j=\int_{C\times\R^+}\car_{\{y\|x\|^\gamma\le P_tK/I\beta_j\}}p_S(y)\, \d y
  \, \d x
\end{equation*}
and $A_0=0$.
\begin{lemma}
  For $j=1,\cdots,\, N_{max}-1$,
  \begin{multline*}
    A_j=\pi R^2Q(\alpha_j-\zeta\ln
    R) \\
    +\pi e^{2/\zeta^2+2\alpha_j/\zeta}Q(\zeta\ln R-2/\zeta-\alpha_j),
  \end{multline*}
where
\begin{equation*}
  \alpha_j=\frac{1}{\sigma}(10\log_{10}(P_tK/I\beta_j)-\mu) \text{ and }
  \zeta=\frac{10 \gamma}{\sigma \ln 10}.
\end{equation*}
\end{lemma}
\begin{proof}
  We can write
  \begin{equation*}
    A_j=\int_C \P(S\|x\|^\gamma \le \tilde{\beta}_j)\, \d x
  \end{equation*}
where $\tilde{\beta}_j=P_tK/I\beta_j$. Remind that $S$ is equal in
distribution to $\exp({\mathcal N}(\mu,\sigma^2)\xi)$ with
$\xi=\ln(10)/10$. Thus after a few manipulations, we get
\begin{equation*}
  A_j=2\pi \int_0^R r \, Q(\alpha_j-\zeta\ln r)\, \d r,
\end{equation*}
where 
\begin{equation*}
  Q(x)=\frac{1}{\sqrt{2\pi}}\int_{-\infty}^x \exp(-\, \dfrac{u^2}{2})\, \d u.
\end{equation*}
The final result follows by a tedious but straightforward integration by parts.
\end{proof}
\begin{thm}
  For any function $\theta\, :\, \R\to \R$,
  \begin{multline*}
    \int \theta(N(x,\ y)) p_S(y)\d
y \ \d \Lambda(x)\\
=\sum_{j=1}^{N_{max}-1}\!\!\theta(j)
(A_j-A_{j-1})+\theta(N_{max})(\pi R^2-A_{N_{max}-1}).
  \end{multline*}
\end{thm}
\begin{proof}
  Since $N$ can take only a finite number of values, we have
  \begin{multline*}
   \int \theta(N(x,\ y)) p_S(y)\d
y \ \d \Lambda(x)\\
=\frac{\rho}{\nu}\sum_{j=1}^{N_{max}}\theta(j)\int_{C\times \R^+}\car_{\{(x,\,
  y),\ N(x,\, y)=j \}} p_S(y)\d
y \ \d x.
  \end{multline*}
Now we see that 
\begin{equation*}
   N(x,\, y)=j \Longleftrightarrow \tilde{\beta}_{j-1}< y
   \|x\|^\gamma\le \tilde{\beta}_j,
\end{equation*}
for $j=1, \cdots, \, N_{max}-1$ and $N(x,\, y)=N_{max}$ when $y
   \|x\|^\gamma >\tilde{\beta}_{N_{max}-1}$. The proof is thus complete.
\end{proof}
We used the  same set of values as for the simulation of Section
\ref{deterministic} together with  assumptions
\ref{a7} and \ref{a8} with $\beta_{min}=0.2$. Results of Table \ref{tab:random}
 show that the theoretical bound  is rather stable when gains become stochastic.

\begin{table}[!ht]
  \centering
  \begin{tabular}{ |c|c|c|c|c|c|c|}
\hline
 $\alpha$  & 1.5 & 1.6 &1.7 & 1.8 & 1.9 & 2\\\hline\hline
$P_{sup}$ &    0.2 &   0.1  &  0.05 &   0.02 & 0.01 & 0.004\\\hline
$\Delta$  & 1.7 & 1.8  & 2.1 & 2.3 & 2.4 & 2.6\\
\hline
  \end{tabular}
   \caption{Comparison between $P_{sup}$ and $P_{loss}$ for
     random gain.}
  \label{tab:random}
\end{table}

\section{General case}
\label{general}
Actually, the method can be applied to more general situations as we
illustrate now.
We  consider only assumptions~\ref{a1}, \ref{a2},
\ref{a3} and~\ref{a7}, a non-frequency selective random gain $G$ with 
distribution $p_G$, a finite number of antennas with a deterministic pattern
and the assumption that the user will receive the signal from the antenna which
can provide a better signal-to-interference ratio. 

Now $C$ is the Borel set where it is possible to find users. Possibly, 
$C=\R^2$. In this region, we have a finite number of antennas $J+1$ 
with the $j$-th located at $y_j$, $j=\overline{1,J}$, and the one we observe 
located at $y_0$. This means that for each user, there is a vector 
$\textbf{G}=(G,G^1,...,G^J)$ where $G$ is the gain due the antenna at $y_0$ and
$G^j$ due to antenna located at $y_j$. We then define the Poisson point process 
$\Psi$ in $C\times \R_+^{J+1}$, representing the user positions and the
gain of each one due to each antenna. Again, since gains from different antennas and positions are
independent altogether, $\Psi$ has intensity $\lambda_{m}$:
\begin{eqnarray*}
\lambda_{m}(\textbf{g},x)=p_{G}(g)p_{G}(g^1)...p_{G}(g^J)\frac{\rho(x)}{\nu}
\end{eqnarray*}
We define the sets
\begin{equation*}
A'=\left\{\bigcup_{j=1}^{J}\left((g,g^1, \ldots,\, x),g^j>\frac{\lVert x-y_0
\rVert^{\gamma}}{\lVert x-y_j\rVert^{\gamma}}g\right)\right\}
\end{equation*}
and
\begin{equation*}
  B=\left((s,x),\,  s\leq R(x)\right),
\end{equation*}
where
\begin{equation*}
  R(x)=\frac{P_tK}{\beta_{min}\lVert x-y_0\rVert^{\gamma}}.
\end{equation*}
The event $((\textbf{G},X)\in A^\prime)$ means that the antenna at $y_0$ provides the highest 
signal-to-interference ratio to a point $X$. The event $((S,X)\in B)$
means the signal-to-interference ratio provided by the antenna
at $y_0$ to a point $X$ is higher than $\beta_{min}$. By Theorem
\ref{thm_concentration}, we are thus led to compute
\begin{equation*}
  \iint_{A\cap B} N(\|x\|, g)^k \d\lambda_m(x,g),
\end{equation*}
for $k=1, \, 2.$ There is no longer a closed form formula for these
integrals but they can be easily and quickly  computed by numerical methods. This
yields  to an upper bound of $P_{loss}$.
We simulated in this section the loss probability for an antenna placed at the
origin and six other antennas placed at the points $y_1=(2R,0)$, 
$y_2=(R,R\sqrt{3})$, $y_3=(-R,R\sqrt{3})$, $y_4=(-2R,0)$, $y_5=(-R,-R\sqrt{3})$
and $y_6=(R,-R\sqrt{3})$, representing an hexagonal arrangement. All other
parameters are the same as the ones in previous sections. We find a mean 
$m_N=21.60$ and the second moment $v_N=26.81$. It turns out that the results are close to the results in Section~\ref{random one},
 suggesting that the approach of
Section~\ref{random one} is satisfactory enough with our physical assumptions.

\section{Concluding remarks}
\label{conc}
Using the concentration and deviation inequalities and the difference 
operator on Poisson space, we have calculated the upper bound probability of 
overloading the system by high demand of subcarriers, over path loss and shadow
fading. To do this we have found the first and second moment of the 
marked Poisson point process of users. 
We conclude that it is possible to find an upper bound for the overloading 
probability, even in a relatively complex system, which is analytically 
computable in a very simple fashion.
The method works for any functional of the configurations, possibly
enriched by marks, which depends only on the positions of each user.
It does not work for functionals involving relative distance
between two or more users. Actually, for such a functional $F$, there is no
bound on $D_xF(\omega)$ valid for all $x$ and $\omega$.

\end{document}